# Webbing and orchestration
## Two interrelated views on digital tools in mathematics education

Luc Trouche, French Institute of Education, ENS de Lyon, France

Paul Drijvers, Freudenthal Institute, Utrecht University, the Netherlands

*Version 23-03-2014*

## 1. Introduction

Since their introduction in the 80s of the previous century, the integration of digital tools in mathematics education is considered both promising and problematic. On the one hand, the potential of digital tools is widely recognized. New technologies are expected to open up new horizons of explorative and productive mathematics education, in which opportunities for experiencing dynamics, for visualizing concepts, and for integrating different representations and resources provide new means of rich and meaningful teaching and learning activities (e.g., see NCTM, 2008; Pea, 1987). On the other hand, the integration of digital tools is considered as problematic in the sense that their availability questions the goals of mathematics education as well as current teaching practices (Lagrange et al., 2003). The latter type of question has not yet been answered in a satisfying way, which may explain the limited integration of digital tools in mathematics education in spite of earlier high expectations.

In order to address the problematic character of the integration of digital tools, educational research was needed that would go beyond enthusiastic reports from first adopters who engaged in small-scale design experiments, and that would be firmly based on theoretical foundations. For some time, such foundations were lacking; in the 90s, however, two important theoretical views emerged: the notion of *webbing* and the *instrumental approach* to tool use. Building on ideas on scaffolding, the notion of webbing was introduced by Noss and Hoyles (1996). Noss and Hoyles were inspired by the work of Papert (e.g., see Papert, 1993) and elaborated his views for the context of using digital learning environments called microworlds for mathematics education. The LOGO programming language served as a paradigmatic environment. In the meantime, the instrumental approach emerged in France (e.g., see Guin & Trouche, 1999). This approach, which has its roots in cognitive ergonomics (Vérillon & Rabardel, 1995) and the anthropological theory of didactics (Chevallard, 1999), highlights the importance of a meaningful relationship between a tool and its user for carrying out a specific task. To stress the importance of the teacher in the students' process of establishing such a relationship, the notion of *instrumental orchestration* was developed (Trouche, 2004).

Initially, the notions of webbing and instrumental orchestration seemed to be rather disconnected, also because of their different cultural and theoretical backgrounds. Whereas the webbing metaphor was oriented towards students' productive and expressive activities, the focus of the instrumental approach was more on the appropriation of given tools for mathematical learning. From the beginning, however, the two 'schools' not only made use of metaphors, but also shared an interest in using digital tools for learning and both highlighted the fact that tools are not just neutral but do impact on student behaviour and on student learning:

> Tools matter: they stand between the user and the phenomenon to be modelled, and shape activity structures. (Hoyles & Noss, 2003, p. 341)

As a consequence, the two views gradually came into contact with each other, and to a certain extent informed each other's further development. It is this mutual influence and development that is at the heart of this article's interest. The central theme of this paper, therefore, is (1) to investigate the ways in which these theoretical notions 'travelled through time', each on its own and in relation to each other, and (2) to assess their impact on future research and teaching. Not claiming a comparison and contrasting that characterises a networking of theories study (e.g., see Prediger, Bikner-Ahsbahs, & Arzarello, 2008; Drijvers, Godino, Font, & Trouche, 2013), we want to describe the distinct and joint journeys of the two perspectives at stake. To do so, the paper follows a more or less chronological line. In section 2, the two 'main actors' are described as they emerged by the end of the previous century. Next, section 3 identifies some differences between the two views, as they became apparent in a symposium in 2003. A further rethinking of the two concepts took place in the frame of the 17$^{th}$ ICMI study (Hoyled and Lagrange, 2010), particularly with respect to the idea of connectivity. As a final development, we focus in section 5 on teachers and teaching, and on the ways the theoretical frameworks inform teachers' professional development. A recent ZDM issue on re-sourcing teacher work, edited by Pepin, Gueudet, and Trouche (2013) serves as a benchmark here.

## 2. The two main actors: webbing and orchestration

In this section we briefly describe the main characteristics of the two main 'actors' in this paper: the notions of webbing and of instrumental orchestration.

### 2.1 Webbing

Noss and Hoyles introduced the notion of webbing in their 1996 book 'Windows on mathematical meanings'. This book was to a large extent inspired by the work of Papert and other on constructionism (e.g., see Papert, 1993 and Papert & Harel, 1991). As a point of departure, the authors see learning as the construction of a web of connections of mathematical ideas, such as connections between classes of problems, between mathematical objects and relationships, or between 'real' entities and personal situation-specific experiences. These connections are captured in intellectual structures, which are individual and can be (re)constructed by students not only through the scaffolding by a teacher but also by themselves through the use of digital environments.

The idea of webbing is a metaphor for this process of building and using such structures. In industry, webbing refers to flexible woven material consisting of a network of natural or artificial fibres, and clearly the reference to flexibility and the interwoven character of the structure fits well Noss and Hoyles' view on learning. The authors themselves describe webbing as follows.

Like the web of mathematical ideas, the Web (we will use a capital to denote the electronic network), is too complex to understand globally – but local connections are relatively accessible. At the same time, one way – perhaps the only way – to gain an overview of the Web is to develop for oneself a local collection of familiar connections, and build from there outwards along lines of one's own interests and obsessions. The idea of webbing is meant to convey the presence of a structure that learners can draw upon and reconstruct for support – in ways that they choose as appropriate for their struggle to construct meaning for some mathematics. (Noss & Hoyles, 1996, p. 108)

In this sense, webbing is considered a fundamental motor for the construction of mathematical meaning, and a step towards situated abstraction. The resulting structures form a support system, a system that the learner has access to and is in control of. Digital environments provide excellent opportunities to build and use such webbing structures, which are shaped by and within the medium. The word 'webbing' refers to both the process of constructing and using connecting structures, and to the resulting structures themselves. As illustrative examples, Noss and Hoyles refer to student observations that they label as flagging, adjusting, sketching, and patterning.

When addressing webbing in the sense of Noss and Hoyles, we should be aware that nowadays, we may have strong associations with the omnipresent Internet, which in the 1990s was only present in its very early and rudimental forms. Some care may be needed to frame the notion of webbing in its current era.

**2.2 Instrumental orchestration**

The instrumental approach to mathematics education emerged in the French context of didacatical engineering, the theory of didactical situations and the anthropological theory of didactics. It initially mainly dealt with complex technological environments such as CASs that were provided to students and teachers, but that were not at all designed for educational purposes (e.g., see Guin & Trouche 1999, Artigue 2002). At its start, it took three essential ideas from cognitive ergonomics (Vérillon & Rabardel, 1995):

1. The distinction between an *artefact* and an *instrument*: an artefact is a material or abstract object, a product of human activity, and is used for performing a type of task. Examples are calculators, and algorithms for solving quadratic equations. These artefacts are *given* to a subject; an instrument is what the subject *builds* from the artefacts;
2. The acknowledgement that the process of appropriating an artefact and building an instrument is a complex one. This process is called instrumental genesis, in which a new entity is given birth. This new instrument is composed of the part of the artefact really mobilized and of a *scheme*, i.e., the structured cognitive organisation of the activity for carrying out the task;
3. The comprehension that instrumental genesis is not a single process, but fundamentally involves two interrelated movements, namely instrumentation (in short: the artefact engraving its mark on the user activity), and instrumentalisation (in short: the user engraving its mark on the artefact).

The proper contribution of the instrumental approach to an understanding of the complexity of tool use in mathematics education was to consider a scheme as a source for building mathematical knowledge. For example, Drijvers (2002, p. 226) describes the emergence of the scheme Isolate-Substitute-Solve for solving a system of two equations ($x + y = 31$ and $x^2 + y^2 = 25^2$) with two unknowns in a symbolic calculator environment (Figure 1). Underlying this scheme is the extending conceptualisation of "what solving an equation is", from "finding a numerical value for an unknown", to "expressing one unknown as a function of another".

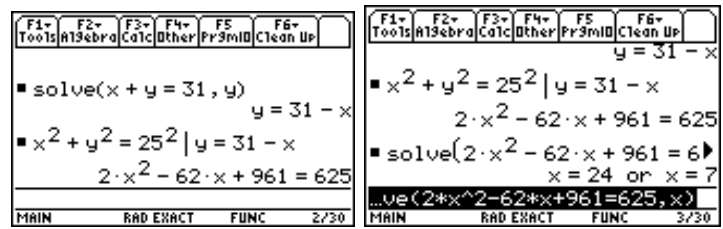

**Figure 1.** *A trace of the scheme Isolate-Substitute-Solve on a TI-89 calculator.*

The meaning of the word "instrumental" in the instrumental approach is quite different from Skemp's (1976) description of the relational – instrumental understanding dialectic. In the latter view, instrumental understanding is "applying rules without reasons". In the instrumental approach, developing techniques is *creating* reasons: instrumentation and conceptualisation are deeply intertwined. This mutual influence leads to taking into account the responsibility of the teacher, and more generally, of the educational institution, for steering students' instrumental geneses. This awareness leads to the notion of instrumental orchestrations, that Trouche (2004, p. 296) defines as follows.

> An instrumental orchestration is defined by *didactic configurations* (i.e., the layout of the artefacts available in the environment, with one layout for each stage of the mathematical treatment) and by *exploitation modes* of these configurations.

The *sherpa-student* configuration depicted in Figure 2 can be considered as paradigmatic of the notion of orchestration. The didactical configuration consists of an organisation of the classroom space, in which each student has his own handheld device, that can be connected to a projection device. All students are seated so that they can see the projection. The exploitation mode comes down to the teacher inviting one student to connect his device and to show his work or to carry out tasks suggested by the teacher. This gives means to the teacher to follow and monitor students' instrumented activity, related to the demands of the mathematical task at stake.

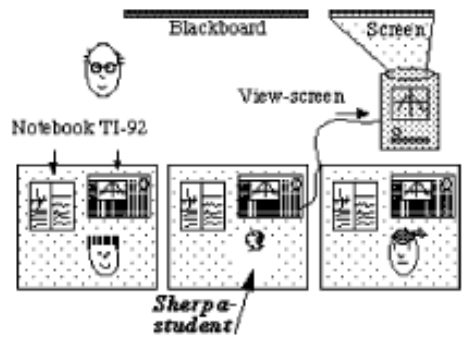

**Figure 2.** *The sherpa-student configuration (Trouche, 2004, p. 298)*

The idea of orchestration constitutes a powerful metaphor that also appears in other theoretical frames. For example, Mariotti and Maracci (2012) include the notion of orchestration in their semiotic mediation approach in the sense of "developing shared meanings, having an explicit formulation, de-contextualized from the artefact use, recognizable and acceptable by the mathematicians' community". The work by Stein and colleagues (2008) focuses at the orchestration of classroom discussions, for "helping teachers move beyond show and tell". In this case, orchestration remains a global perspective and the authors mainly highlight the general conditions for such discussions to happen: they make a plea for designing learning environments "so that students are "authorized" to solve mathematical problems for themselves, are publicly credited as the "authors" of their ideas, and develop into local "authorities" in the discipline" (Stein *et al.*, 2008, p. 332). Dillenbourg and Jermann (2010) also use the metaphor of orchestration, and propose a musical notation for the sequence of different configurations (individual, group, class, et cetera) aiming to foster student activity (Figure 3). Hähkiöniemi (2013) suggests that teachers need support to develop an appropriate repertoire of orchestrations.

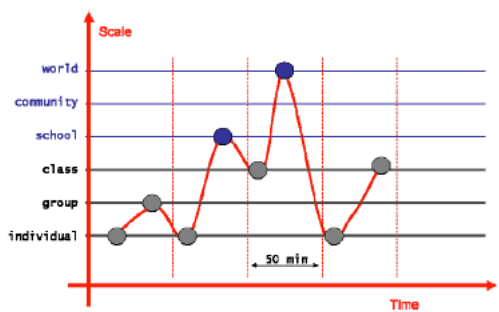

**Figure 3.** *A musical notation for integrated scripts (Dillenbourg & Jermann, 2010)*

The notion of instrumental orchestration as proposed by Trouche (2004, 2005) in mathematics education applies not only to the episodes of discussions and seems to be more precise, aiming to define the management of both space and time of a learning environment, according to the different stages of the task to be carried out.

Webbing as a model of internal students' cognitive functioning on the one hand, and instrumental orchestration as a model of guiding students' instrumental genesis on the other, seem to be at two different levels of mathematics education interests. However, these two concepts have come to meet and to engage dialogue on the occasion of the third CAME[1] symposium, which will be discussed in the next section.

## 3. Differentiating the concepts of webbing and orchestration

The CAME symposium, and the theme Mind and Machine in particular, constituted a fruitful context for discussing the concepts of webbing and orchestration. In a presentation by Trouche and a reaction by Hoyles[2] addressed the approaches bearing respectively the concepts of webbing and orchestration. These two contributions share essential starting points, such as (1) recognizing the "underestimation of the complexity of instrumentation processes" (Hoyles *et al.*, 2004, p. 311), (2) evidencing "a failure to theorize adequately the complexity of supporting learners to develop a fluent and effective relationship with technology in the classroom" (ibidem, p. 311), and (3) underlining the fact that "the marginalisation of technology by educational institutions has also turned our attention to the need for a more precise analysis of the role of the teacher in these new and changing didactical contexts" (ibidem, p. 313). This being said, the two papers reveal three important differentiating nuances.

### 3.1 The true nature of mathematical conceptualisation

In the webbing approach, conceptualisation appears as a *coordination process*, "the process by which the student infers meaning by *coordinating* the structure of the learning system (including the knowledge to be learned, the learning resources available, prior student knowledge and experience and constructing their own scaffolds by interaction and feedback)" (Hoyles, Noss, & Kent, 2004, p. 319).

In the instrumental orchestration approach, conceptualisation appears as a *command process*, characterized by the conscious attitude to consider, with sufficient objectivity, all the information immediately available not only from the calculator, but also from other sources and to seek mathematical consistency between them (Guin & Trouche, 1999). "Very sophisticated artefacts such as the artefacts available in a computerized learning environment give birth to a set of instruments. The articulation of this set demands from the subject a strong *command process*. One of the key elements for a successful integration of these artefacts into a learning environment is the institutional and social assistance to this individual command process. Instrumental orchestrations constitute an answer to this necessity." (Trouche, 2004, p. 304). It seems that there is a kind of intended *internalization* from an *instrumental orchestration*, seen as an external process of monitoring students' instruments by the teacher, to an *internal orchestration*[3], seen as a process of self-monitoring the individual and personal instruments by a student.

---

[1] The Third CAME (Computer Algebra in Mathematics Education) symposium was held in Reims (France) in 2003, focusing on "Learning in a CAS Environment: Mind-Machine Interaction, Curriculum & Assessment".
[2] Contributions are available on the conference website (http://www.lkl.ac.uk/research/came/events/reims/). The papers by Hoyles, Noss, and Kent (2004) and Trouche (2004) emerged from these contributions.
[3] Using the same metaphor, Dehaene (2010, p. 221) identifies the internal orchestration of distributed neuronal networks as a central question for the neurosciences.

Coordination and control are certainly two facets of mathematical activity, particularly in technological rich environment, and the two approaches seem to privilege, each, one of these facets.

**3.2 The role of interactions between learners**

Hoyles *et al.* clearly situate this nuance: "Perhaps the different unfolding histories of our respective research efforts – including the choice of technologies to study, has resulted in differences in emphasis between us rather than in core perspective, most notably regarding the weight we accord to the role of interactions among learners [...] (ibidem, p. 317)". Here seems to emerge a kind of intended *externalization*, from the webbing as an internal process (§ 1.1) to a webbing as connections between learners, with mathematical knowledge emerging from these connections. In this approach, social aspects seem to relate mainly to interactions between pupils, while, for the instrumental approach, social aspects are in many cases orchestrated by the teacher. For example, the sherpa-student orchestration depicted in Figure 2 clearly has a social aspect: students will interact on the work shown by the sherpa student under the guidance of the teacher.

Things are actually more complex. For example, in the mirror configuration (Figure 4), a pair of students learns from observing another pair of students solving a problem, without the teacher's intervention, but within an orchestration that was designed by her.

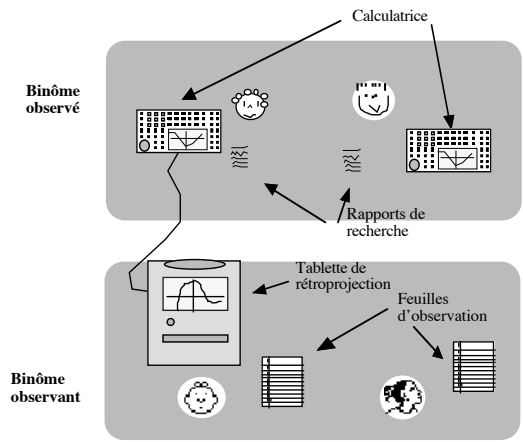

**Figure 4.** *The mirror configuration (Trouche, 2005)*

Overall, there clearly is a difference between Hoyles's approach, speaking of *learners*, interacting in the frame of *communities of practice*, and Trouche's approach, speaking of *students*, interacting in the frame of classroom didactical configurations designed by a teacher. Hoyles *et al*. (ibidem, p. 316) outline this as follows.

> Orchestration/instrumentation is a felicitous metaphorical pairing, which captures something of the relationships involved. One striking issue about the metaphor (although bearing in mind that a metaphor should not be interpreted too literally) is this: what is the analogue of music in the orchestration of mathematical instruments? Presumably the answer is mathematics, but what kind of mathematics? The production of ''official'' mathematics which the teacher (i.e., conductor) wishes ''to hear'' in order to judge the students' performance, or the more personal and community-based norms tied to students' own developing conceptions and personal inclinations?

"Who orchestrates" is certainly an essential question. It seems to constitute a third terrain for confronting the two approaches.

**3.3 The role of the teacher in orchestration**

For the instrumental orchestration, the teacher has a central role, as "the orchestra conductor rather than a one-man's band" (Trouche, 2004, p. 299). The orchestra metaphor allows some variation, playing on the kind of orchestra – the jazz band offers room for student improvisation –, the choice of the score, the choice of the didactical configuration. But in each of these cases, the teacher is always central.

From the webbing perspective, the teacher can remain distant, and digital technology in itself seems to have the potential to orchestrate student activity, as phrased by (Hoyles *et al.*, 2004, p. 319):

> We simply wish to raise the profile of the potential role of the technological tools and the kinds of symbolic language and interactions that they might foster among learners. Consider, just as a simple example, a learning environment that typically develops where rows of computers are physically fixed one behind the other with little if any connectivity, and contrast this with one where students have wireless-connected laptops and can enjoy a hugely greater freedom to collaborate and share their ideas, both face-to-face and virtually. This process of breathing life into technology necessarily differs from one individual to another, if only because of differences between the kinds of things a student takes for granted, already knows, or is trying to understand. In so far as the artefact's properties can be thought of as ''orchestrating'' the actions and expressions of an individual or a group of learners, this orchestration is not invariant across different individuals.

The notions of webbing and orchestration aim to describe the way in which learners may develop a fluent and effective relationship with mathematics through technology, orchestration adding the role of teacher to support such a development.

**3.4 A possible unified point of view**

The two perspectives of webbing and orchestration seem to be engaged in a process of convergence, as Hoyles *et al.* (ibidem, p. 322) acknowledge:

> Our two perspectives are united in recognising that if orchestration aims merely to bring about ''convergence'' of mathematical expression with official mathematical discourse, the potential of the technology will almost certainly be missed. Yet convergence is important, so two qualifications of the previous comment are in order. First, convergence for the students may take time – significant time over years rather than days or months. Second, the process of orchestration must take place at different levels, separated in time: the first level of orchestration being to foster the growth of situated abstractions (or, instrumented social and individual mathematical schemes) which establish a ''cognitive scaffolding'' for a second level of orchestration to bring about convergence or at least alignment through discussion of boundary objects. This second phase might be expected to take place over an extended period, and through a combination of collective activity in the classroom and individual work by students. In fact, both notions highlight a view on learning that considers mental constructions, i.e. schemes or structures, as crucial, and both see a potential for digital technology to support the construction process, by means of scaffolding and processes such as instrumentation and instrumentalization.

The CAME conference discussion took place within a "circle of experts", but had an impact on a larger community of research interested in digital technology in mathematics education (Lagrange *et al.*, 2003). This led, at an international level, to "rethinking the terrain".

# 4. Rethinking the terrain… and the concepts

The seventeenth ICMI Study entitled "Technology Revisited" had its study conference in Hanoï in 2006[4]. The title of its proceedings "Mathematics Education and Technology—Rethinking the Terrain" (Hoyles & Lagrange, 2010) reveals the deep evolution of this research domain from the first naive approach towards an acknowledgment of the complexity of integrating digital technology in regular mathematics classrooms, that also affects the notions of webbing and orchestration. This study opened new horizons related to the impact of Internet on learning and teaching: the opportunities for collaborative work between students, but also between teachers, thus evidencing the crucial role of *connectivity*. We will look at these evolutions through the window of the panel on connectivity and virtual networks held during this conference, which was chaired by Celia Hoyles (Hoyles *et al.,* 2010). This panel proposes five contributions, including one by Trouche and Hivon entitled "Connectivity: new challenges for the ideas of webbing and orchestrations".

### 4.1 Orchestrations as living entities

Trouche and Hivon's contribution describes the collaborative work of six teachers who carefully designed an orchestration of a mathematical activity in a calculator network environment. The organisation of tables and calculators (Figure 5) as well as a relevant choice of the various applications available are discussed.

---

[4] Its Discussion Document « Technology revisited » is available at
http://www.mathunion.org/fileadmin/ICMI/files/Digital_Library/DiscussionDocs/DD_icmiStudy17_02.pdf.

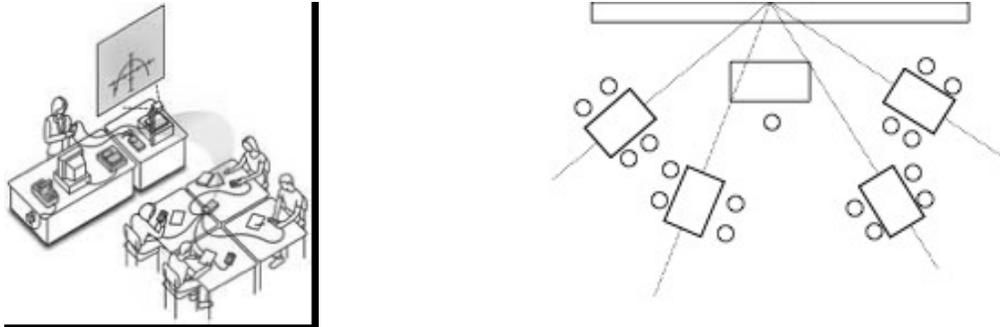

**Figure 5.** *From the intended configuration by the manufacturer to the radius configuration chosen by the teachers*

Orchestrations appear here as *living entities* in two respects:

- In the teachers' respect, orchestrations appear as the result of their collaborative work of together designing a priori orchestrations and of analysing the classroom implementations: "The cross-observations of teachers in their own classes helped them to create a distance from their own practice and to develop a reflexive attitude to the orchestration of students' activity" (Hoyles *et al.*, 2010, p. 449);
- In the students' respect, students no longer just play the music written by the conductor; rather they are composing part of the music themselves. The question, then, is how the teacher can create conditions to make the music not too different from his/her intentions, and to enrich the prepared partition with the – sometimes unexpected – student improvisations.

As Dillenbourg and Jermann (2010, p. 527) claim, "The key difference between music orchestration and classroom orchestrations is that, when orchestrating a classroom, the score has often to be modified on the fly". As the initial notion of orchestration included didactical configurations and exploitations modes (§ 2.2), and as such did not focus on the dynamic adaptions made *on the fly*, Drijvers, Doorman and colleagues (2010) introduced a third orchestration level of *didactical configuration.* If we see orchestrations as dynamic entities that emerge in cycles of preparation, adaptation on the fly and reflection, there is an interesting parallel between students engaging in instrumental genesis and teachers engaging in process of developing a repertoire of instrumental orchestration; a process of orchestrational genesis, so to say.

**4.2 Orchestrating to foster situated abstractions and negotiate shared meanings**

Comparing the sherpa configuration (Figure 2) and the radius configuration (Figure 5), the connectivity in the latter one clearly affords more flexible relationships between pupils. For Trouche and Hivon (Hoyles et al., p. 448), through connectivity, "a new interactivity is fostered between the artefact and the student, and between students themselves: students convey their messages through the artefact, the artefact acts on the students enabling them to distance themselves from their productions thus freeing them to become more easily involved in peer exchanges. Thus the common space (Figure 6) becomes a space of debate and exchange that aims to elaborate a social 'mathematical truth'."

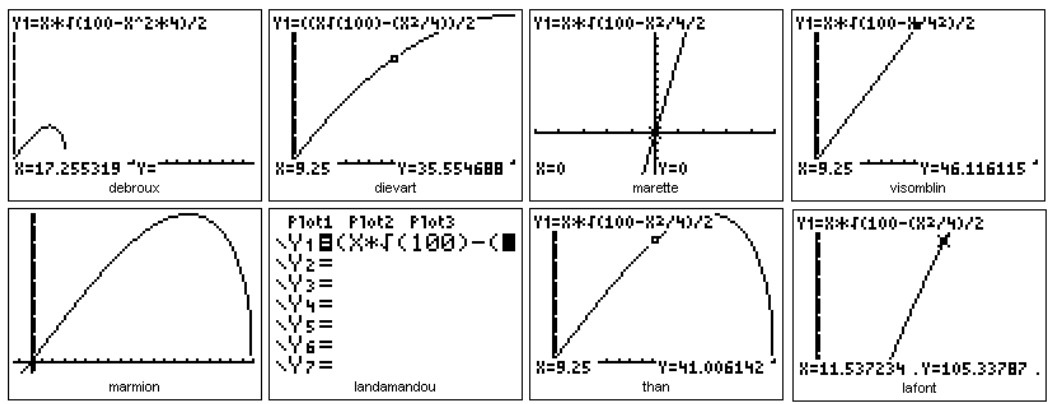

**Figure 6.** *The classroom screen displaying the different student calculator screens becoming a common work space*

Such a connected configuration invites a flexible view on mathematics learning, in which shared meanings (in this case on the notion of function) emerge. Finally, Trouche and Hivon (ibidem, p. 451) stress "the need to rethink the notion of orchestration and the notion of webbing […]. On the one hand, it is important to have in mind a necessary assistance (the notion of orchestration) of students' mathematical activity, and on the other, it is crucial to let the students free to think and establish connections (the idea of webbing)".

The ICMI Study panel suggests some converging points of view, both on the true nature of mathematical conceptualisation (§ 3.1) and on the power of interactions between learners (§ 3.2). The importance of the teacher (§ 3.3) is largely acknowledged by Hoyles (ibidem, p. 460): "[…] the teacher is crucial. But here we are delineating new, even more demanding roles for the teacher, to be aware – across not only her own classroom but those in remote locations – of the evolution of discussion, the mathematical substance of what is and what is not discussed, and the need all the while to find ways to keep students on task without removing the exploratory and fun elements of the work". These demanding roles for the teacher suggest an urgent research focus, in which teachers' development is considered: "Teachers are themselves involved in a process of instrumental genesis to develop artifacts into instruments for accomplishing their teaching tasks" (Drijvers, Kieran, & Mariotti, 2010, p. 112). The re-sourcing of teacher work is the topic of the next section.

## 5. RE-SOURCING TEACHER WORK

In the previous sections, we noticed that it is naïve to expect the process of webbing, which shares characteristics with the instrumental genesis central in the instrumental approach, to be autonomous and self-evident. Rather, it needs guidance by the teachers, and this is where the notion of instrumental orchestration fits in. To be able to provide this guidance, however, teachers' professional development is needed. The question, therefore, is how we can put into practice the notions of webbing and orchestration to study teacher behaviour and to inform teachers' professional development. This questionis addressed in the ZDM issue 45(7) entitled 'Re-sourcing teacher work and interaction: new perspectives on resource design, use and teacher collaboration', edited by Pepin, Gueudet and Trouche (2013). Three contributions to this benchmark issue will be discussed in below.

## 5.1 Towards a taxonomy of orchestrations?

The notion of instrumental orchestration can be used for different purposes. Whereas Trouche's (2004) initial paper on instrumental orchestration highlights the potential and the importance of the concept, Drijvers used instrumental orchestration primarily to observe teaching practices put into action by novice and expert teachers and to establish connections between these practices and the teachers' views on the role of digital tools in mathematics education (Drijvers, Doorman, Boon, Reed, & Gravemeijer, 2010). In a paper in the ZDM issue discussed here, this led to a tentative taxonomy of teaching practices, which is not claiming completeness; rather, a different teaching setting with other digital tools might lead to a different inventory (Drijvers, Tacoma *et al.*, 2013; see Figure 7). A small-scale follow-up study suggests that differences between expert and novice teachers (with respect to their experience of integrating of digital technology in teaching) seem to lie more in the quality and duration of the orchestrations than in the richness of the repertoire, even if novice teachers tend to be careful with using new, more 'adventurous' orchestrations in which they might feel less in control (Kaper, 2013, unpublished master thesis). For example, expert teachers were more efficient in dealing with technology-centred orchestrations, such as Technical-demo and Technical-support (see Figure 7). This highlights the importance of the quality of the orchestration, which so far remains largely unexplored.

Ruthven (2014) describes these different purposes, as he perceives them among the co-authors of the present paper, as follows:

> This shift in the meaning and structure of "instrumental orchestration" between Trouche and Drijvers reflects their differing purposes: while Trouche is interested in examining potential strategies, Drijvers is seeking to describe observed patterns. (Ruthven, 2014).

The Drijvers et al. paper not only shows a way of identifying teaching strategies with digital technology; it also shows how a tentative way of putting into practice a theoretical notion may impact on the concept itself.

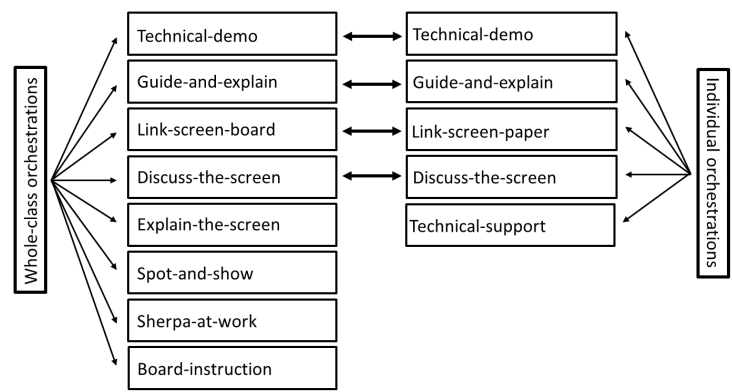

**Figure 7.** *Overview of observed orchestrations (Drijvers, Tacoma et al., 2013)*

## 5.2 The Cornerstone Mathematics Project

In the same ZDM issue, Hoyles and colleagues report on the Cornerstone Mathematics Project, which addresses teachers' professional development through their participation in the project (Hoyles, Noss, Vahey, & Roschelle, 2013). In particular, teachers' adaptation of an existing technology-rich teaching sequence, their ownership and the impact of the collaboration with colleagues are investigated. It is interesting to consider this article's theoretical and methodological approach and its results, and to relate these to the notions of webbing and orchestration, and to other work presented in the same journal issue.

Concerning the theoretical orientation, it is remarkable that notion of webbing is absent in the paper, even if one of the foci is students' conceptual development, as a result of using digital technology. Interestingly, the authors do mention instrumental genesis and describe instruments as "systems with which the user gains fluency and expressive competence" (p. 1058). We wonder if the word 'system' in the previous quote could be related to the word 'structure' in description of webbing in the Noss and Hoyles 1996 book: "the structure that learners can draw upon and reconstruct for support"? In the 2013 paper we discuss here, the notions of orchestration and instrumentalisation are used, even if their meanings may in some cases be quite specific.

In addition to similarities in the theoretical framework, the Hoyles *et al*. 2013 study also makes methodological choices that are shared by others. For example, compared to the Drijvers *et al*. study reported in the same journal volume (Drijvers, Tacoma, Besamusca, Doorman, & Boon, 2013), the two studies share the use of an existing teaching sequence as the point of departure for a new sequence and a corresponding process of adaptation and professional development. Also, small-size communities of teachers are set up, and the communication within these communities is organized in a blended way through both face-to-face meetings and online facilities.

Finally, there is also some agreement on the findings reported in Hoyles et al. (2013) and in Drijvers, Tacoma et al. (2013). For example, it appeared difficult to make teachers engage in online communication, in addition to face-to-face meetings. Also, the importance of ownership that teachers should feel with respect to the teaching unit and its adaptation is a common conclusion. Adaptation does not necessarily imply adoption. Finally, the role of the community of colleagues is highlighted. This relates to the topic of teachers' documentational work and collaborative work, which will be addressed in the next section.

## 5.3 The documentational and collaborative work of teachers

A third contribution in the ZDM issue relates to teachers' collective work with resources (Gueudet, Pepin, & Trouche, 2013). The point of departure is that teachers nowadays have access to an immense and ever growing collection of teaching resources. Even if the internet is an important vehicle for dissemination, these resources are not necessarily technology related: interactive applets can be complemented by paper-on-screen documents such as book chapters or paper-and-pencil tasks. This myriad of resources now confronts the teacher with the challenge of finding, selecting, adapting and using the available teaching resources for the purpose of their own teaching. Is digital tools are provided to students, teaching resources are provided to teachers, and it is here that an analogy with instrumental genesis emerges: the process of teachers transforming resources into documents for their own teaching is called documentational genesis (Gueudet & Trouche, 2009).

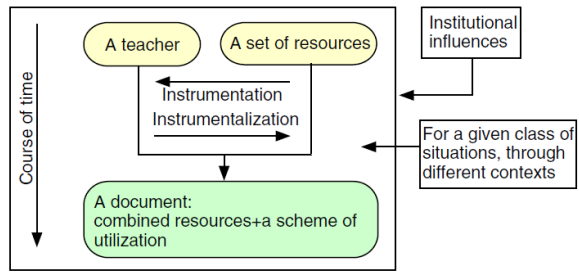

**Figure 8.** *A representation of documentational genesis (Gueudet & Trouche, 2009)*

This process, which teachers engage in, is depicted in Figure 8. Sabra (2011) exemplified the kind of resource systems teachers can develop during their documentational genesis (Figure 9). Kieran and colleagues (2013) recently applied the notion of documentational genesis also to the work of the educational researcher.

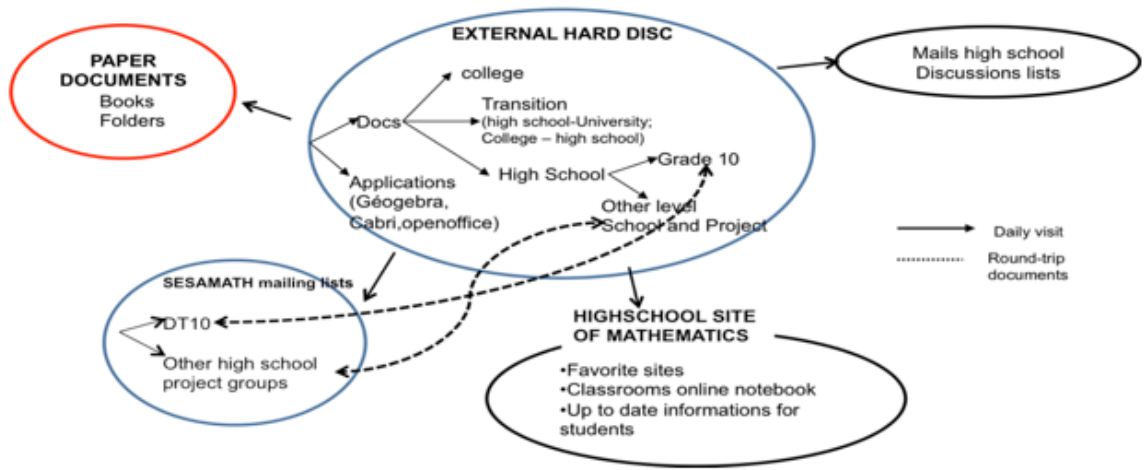

**Figure 9.** *A self representation of one teacher's system of resources (Sabra & Trouche, 2011)*

It should be noticed that, in case digital tools are involved, documentational genesis includes the teacher's personal instrumental genesis as well as her orchestrational genesis (see 4.1): a teachers' documentational activity includes a grounding instrumental genesis concerning the artifacts in use, and is followed by the design of instrumental orchestrations aiming at the integration of those artefacts in the teaching processes (grounding the orchestrational geneses). Finally, this process ends with the documentational geneses, fed by all the resources that the teacher encountered and gathered for her own needs. As is the case for instrumental genesis, these documentational geneses are complex and not easy to pursue in isolation. Therefore, it seems that collaboration among teachers, in which orchestrations are shared and co-designed, can be very fruitful. It is here that the notion of community of practice comes in. Even if it is not self-evident that groups of teachers working together form an efficient community of practice (Besamusca & Drijvers, 2013), the initiatives Hoyles and colleagues undertook to establish collaborative work by teacher through the NCETM and in the Cornerstone Mathematics project (Hoyles et al., 2013) are promising.

As a conclusion of this section, the three papers mentioned here suggest some convergence, both from a theoretical and a methodological perspective, as well as in results. As an overall lesson, we learn that in-service teachers' professional development can be fostered by a joint work on technology-rich teaching sequences for which teachers feel ownership, and which is done in small-scale communities through blended means of communication.

# 6. FUTURE DIRECTIONS TOWARDS EXPLORING TERRA INCOGNITA

How do the above reflections on the notions of webbing and orchestration guide future research agendas? Teachers' collaborative work seems to constitute a true terra incognita to be explored: "The way digital technologies can support and foster today collaborative work, at a distance or not, between students or between teachers, and also between teachers and researchers, and the consequences that this can have on students' learning processes, on the evolution of teachers' practices is certainly one essential evolution that educational research has to explore in the future. As mentioned above, most of this space is still for us nearly terra incognita" (Artigue, in Hoyles & Lagrange, 2010, p. 473).

Collective work with resources appears to be "an essential dimension for teacher documentation" (Gueudet, Pepin, & Trouche, 2013), and for this reason an essential field for future research. The emergence of online teacher associations, such as Sésamath in France (Sabra & Trouche, 2011) provides interesting opportunities to analyse and better understand the interaction between individual and collective documentational geneses.

New contexts for learning are offered by MOOCs (Massive Open Online Courses), which confront teams of teachers to, in many cases, several thousands of students. These new formats allow students to learn through different types of complex social interactions on the web, and, as such, challenge the usual repertoire of orchestrations. This development forms a new opportunity to rethink the relationships between orchestrations and webbing.

For investigating these new fields, new methodologies are needed, as well as new theories: "…teachers's actions in supporting new communities of practice are recognised as crucial, and new roles for the teacher noted, although it is acknowledged that these roles have as yet been undertheorised" (Hoyles & Lagrange, 2010, p. 423). It is this fascinating field of technological, social and educational developments that we explore, gratefully capitalizing on the work done by Celia Hoyles and her colleagues over the last decades.

# LITERATURE

Artigue, M. (2002). Learning Mathematics in a CAS Environment: The Genesis of a Reflection about Instrumentation and the Dialectics between Technical and Conceptual Work. *International Journal for Computers in Mathematical Learning*, *7*(3), 245-274.

Chevallard, Y. (1999). L'analyse des pratiques enseignantes en théorie anthropologique du didactique. *Recherches en Didactique des Mathématiques, 19*(2), 221-266.

Dehane, S. (2010). *La bosse des math : 15 ans après*. Paris : Odile Jacob.

Dillenbourg, P., & Jermann, P. (2010). Technology for classroom orchestration. In M. S. Khine, & I. M. Saleh (Eds.), *New Science of Learning: Cognition, Computers and Collaboration in Education* (pp. 525-552). New York / Berlin: Springer.

Drijvers, P. (2002). L'algèbre sur l'écran, sur le papier et la pensée algébrique. In D. Guin, & L. Trouche (Eds.), *Calculatrices symboliques. Faire d'un outil un instrument du travail mathématique : un problème didactique* (pp. 215-242). Grenoble : Editions La Pensée sauvage*.*

Drijvers, P., Doorman, M., Boon, P., Reed, H., & Gravemeijer, K. (2010). The teacher and the tool: instrumental orchestrations in the technology-rich mathematics classroom. *Educational Studies in Mathematics, 75*(2), 213-234.

Drijvers, P., Godino, J.D., Font, V., & Trouche, L. (2013). One episode, two lenses; A reflective analysis of student learning with computer algebra from instrumental and onto-semiotic perspectives. *Educational Studies in Mathematics, 82*(1), 23-49.

Drijvers, P., Kieran, C., & Mariotti, M. A. (2010). Integrating technology into mathematics education: theoretical perspectives. In C. Hoyles, & J.-B. Lagrange (Eds.), *Mathematics education and technology - Rethinking the terrain* (pp. 89-132). New York/Berlin: Springer.

Drijvers, P., Tacoma, S., Besamusca, A., Doorman, M., & Boon, P. (2013). Digital resources inviting changes in mid-adopting teachers' practices and orchestrations. *ZDM, The International Journal on Mathematics Education, 45*(7), 987-1001.

Gueudet, G., & Trouche, L. (2009). Towards new documentation systems for mathematics teachers? *Educational Studies in Mathematics, 71*, 199–218.

Gueudet, G., Pepin, B., & Trouche, L. (2013). Collective work with resources: an essential dimension for teacher documentation, *ZDM, The International Journal on Mathematics Education, 45*(7), 1003-1016.

Guin, D., Ruthven, K., & Trouche, L. (Eds.) (2005). *The didactical challenge of symbolic calculators: turning a computational device into a mathematical instrument.* New York/Berlin: Springer.

Guin, D., & Trouche, L. (1999). The complex process of converting tools into mathematical Instruments. The case of calculators. *International Journal of Computers for Mathematical Learning, 3*, 195–227.

Hähkiöniemi, M. (2013). Teacher's reflections on experimenting with technology-enriched inquiry-based mathematics teaching with a preplanned teaching unit. *Journal of Mathematical Behavior, 32*, 295-308.

Hoyles, C., Kalas, I., Trouche, L., Hivon, L., Noss, R., & Wilensky, U. (2010). Connectivity and Virtual Networks for Learning. In C. Hoyles, & J.-B. Lagrange (Eds.), *Mathematical Education and Digital Technologies: Rethinking the terrain* (pp. 439-462)*.* New York/Berlin: Springer.

Hoyles, C., & Lagrange, J.B. (Eds.) (2010). *Mathematical Education and Digital Technologies: Rethinking*

*the terrain.* New York / Berlin: Springer.

Hoyles, C., & Noss, R. (2003). What can digital technologies take from and bring to research in mathematics education? In A.J. Bishop, M.A. Clements, C. Keitel, J. Kilpatrick, & F. Leung (Eds.), *Second international handbook of mathematics education, Vol. 1* (pp. 323-349). Dordrecht, the Netherlands: Kluwer Academic Publishers.

Hoyles, C., Noss, R., & Kent, P. (2004). On the integration of digital technologies into mathematics classrooms. *International Journal of Computers for Mathematical Learning, 9*, 309–326.

Hoyles, C., Noss, R., Vahey, P., & Roschelle, J. (2013). Cornerstone Mathematics: designing digital technology for teacher adaptation and scaling. *ZDM, The International Journal on Mathematics Education, 45*(7), 1057-1070.

Kieran, C., Boileau, A., Tanguay, D., & Drijvers, P. (2013). Design researchers' documentational genesis in a study on equivalence of algebraic expressions. *ZDM, The International Journal on Mathematics Education, 45*(7), 1045-1056.

Lagrange, J.-B., Artigue, M., Laborde, C., & Trouche, L. (2003). Technology and Mathematics Education: a Multidimensional Study of the Evolution of Research and Innovation. In A.J. Bishop, M.A. Clements, C. Keitel, J. Kilpatrick, & F.K.S. Leung (Eds.), *Second International Handbook of Mathematics Education* (pp. 239-271). Dordercht, the Netherlands: Kluwer Academic Publishers.

Mariotti, M.A., & Maracci, M. (2012). Resources for the teacher from a semiotic mediation perspective. In G. Gueudet, B. Pepin, &, L. Trouche (Eds), *From text to 'lived resources': curriculum material and mathematics teacher development* (pp. 59-75). New York: Springer

National Council of Teachers of Mathematics (2008). *The role of technology in the teaching and learning of mathematics.* http://www.nctm.org/about/content.aspx?id=14233. Accessed 31 July 2013.

Noss R., & Hoyles C. (1996). *Windows on Mathematical Meanings. Learning Cultures and Computers*. Dordercht, the Netherlands: Kluwer Academic Publishers.

Papert, S. (1993). *Mindstorms: Children, Computers and Powerful Ideas*. New York : Basic Books.

Papert, S. & Harel, I. (1991). *Situating Constructionism*.
Retrieved from http://www.papert.org/articles/SituatingConstructionism.html.

Pea, R. (1987). Cognitive technologies for mathematics education. In A.H. Schoenfeld (Ed.), *Cognitive Science and Mathematics Education* (pp. 89-122). Hillsdale: Lawrence Erlbaum .

Pepin, B., Gueudet, G., & Trouche, L. (Eds.) (2013). Re-sourcing teacher work and interaction: new perspectives on resource design, use and teacher collaboration, *ZDM, The International Journal on Mathematics Education, 45*(7).

Prediger, S., Bikner-Ahsbahs, A., & Arzarello, F. (2008). Networking strategies and methods for connecting theoretical approaches: first steps towards a conceptual framework. *ZDM, The International Journal on Mathematics Education, 40*(2), 165-178.

Ruthven, K. (2002). Instrumenting mathematical activity: reflections on key studies of the educational use of computer algebra systems. *International Journal of Computers for Mathematical Learning*, *7*(3), 275–291.

Ruthven, K. (2013). From design-based research to re-sourcing 'in the wild': reflections on studies of the co-evolution of mathematics teaching resources and practices. *ZDM, The International Journal on Mathematics Education, 45*(7), 1071-1079.

Ruthven, K. (2014). Frameworks for analysing the expertise that underpins successful integration of digital technologies into everyday teaching practice. In A. Clark-Wilson, O. Robutti, & N. Sinclair (Eds.) *The Mathematics Teacher in the Digital Era* (pp. 373-393). New York/Berlin: Springer.

Sabra, H., & Trouche, L. (2011). Collective design of an online math textbook: when individual and collective documentation works meet. In M. Pytlak, T. Rowland, & E. Swoboda (Eds.)*, Proceedings of CERME 7* (pp. 2356-2366) *,* 9th to 13th February 2011. Rzesów, Poland

Skemp, R. (1976). Relational Understanding and Instrumental Understanding. *Mathematics Teaching*, *77*, 20–26.

Trouche, L. (2004). Managing the Complexity of Human :Machine Interactions in Computerized Learning Environments: Guiding Students' Command process Through Instrumental Orchestrations. *International Journal of Computers for Mathematical Learning*, *9*, 281–307.

Trouche, L. (2005). Construction et conduite des instruments dans les apprentissages mathématiques : nécessité des orchestrations. *Recherches en didactique des mathématiques*, *25*, 91-138.

Trouche, L., & Drijvers, P. (2010). Handheld technology for mathematics education, flashback to the future. *ZDM, The International Journal on Mathematics Education, 42*(7), 667-681.

Vérillon, P., & Rabardel P. (1995). Cognition and artefacts: A contribution to the study of thought in relation to instrument activity. *European Journal of Psychology in Education, 9*(3), 77-101.